\documentclass[12pt]{amsart}
\usepackage{amscd,amssymb}
\usepackage[all]{xy}
\usepackage{amsmath}
\usepackage{graphicx} 

\theoremstyle{plain}
\newtheorem{thm}[subsection]{Theorem}
\newtheorem{lem}[subsection]{Lemma}
\newtheorem{prop}[subsection]{Proposition}
\newtheorem{cor}[subsection]{Corollary}

\theoremstyle{definition}
\newtheorem{rk}[subsection]{Remark}
\newtheorem{definition}[subsection]{Definition}
\newtheorem{ex}[subsection]{Example}

\numberwithin{equation}{section}
\newcommand{\OO}{{\mathcal O}}
\newcommand{\D}{{\mathcal D}}
\newcommand{\cH}{{\mathcal H}}
\newcommand{\I}{{\mathcal I}}
\newcommand{\J}{{\mathcal J}}
\newcommand{\ZZ}{{\mathcal Z}}
\newcommand{\T}{{\mathcal T}}

\newcommand{\M}{{\mathcal M}}
\newcommand{\F}{{\mathcal F}}
\newcommand{\A}{{\mathcal A}}
\newcommand{\X}{{\mathcal X}}

\newcommand{\B}{{\mathcal B}}
\newcommand{\CC}{{\mathcal C}}

\renewcommand{\SS}{{\mathcal S}}
\newcommand{\Specan}{{{\mathcal S}pecan}}

\newcommand{\Z}{\mathbb{Z}}
\newcommand{\Q}{\mathbb{Q}}
\newcommand{\R}{\mathbb{R}}
\newcommand{\C}{\mathbb{C}}
\newcommand{\N}{\mathbb{N}}

\newcommand{\bL}{\mathbb{L}}

\newcommand{\tH}{\widetilde{H}}

\DeclareMathOperator{\rank}{rank}

\DeclareMathOperator{\codim}{codim}
\DeclareMathOperator{\supp}{supp}

\DeclareMathOperator{\can}{can}
\DeclareMathOperator{\Sp}{Sp}

\begin{document}

\title [Multiplier Ideals, V-filtrations and Transversal Sections ]
{ Multiplier Ideals, V-filtrations and Transversal Sections    }

\author[A. Dimca]{A. Dimca }
\address{  Laboratoire J.A. Dieudonn\'e, UMR du CNRS 6621,
                 Universit\'e de Nice-Sophia Antipolis,
                 Parc Valrose,
                 06108 Nice Cedex 02,
                 FRANCE.}
\email
{dimca@math.unice.fr}

\author[Ph. Maisonobe]{ Ph. Maisonobe }
\address{  Laboratoire J.A. Dieudonn\'e, UMR du CNRS 6621,
                 Universit\'e de Nice-Sophia Antipolis,
                 Parc Valrose,
                 06108 Nice Cedex 02,
                 FRANCE.}
\email
{phm@math.unice.fr}

\author[M. Saito]{M. Saito }
\address{ RIMS Kyoto University, Kyoto 606--8502 JAPAN}
\email
{ msaito@kurims.kyoto-u.ac.jp}

\author[T. Torrelli]{T. Torrelli}
\address{  Laboratoire J.A. Dieudonn\'e, UMR du CNRS 6621,
                 Universit\'e de Nice-Sophia Antipolis,
                 Parc Valrose,
                 06108 Nice Cedex 02,
                 FRANCE.}
\email
{torrelli@math.unice.fr}



\begin{abstract} We show that the restriction to a smooth transversal section
commutes to the computation of multiplier ideals and V-filtrations. As an application we prove the constancy of the spectrum along any stratum of a Whitney regular stratification.

\end{abstract}

\maketitle

\section{Introduction} \label{s0}

Let $X$ be a complex $n$-dimensional manifold and $D \subset X$ an effective divisor defined by a holomorphic function $f$. To put the results of this paper in their proper perspective, we start by recalling the local topological triviality of Whitney regular stratifications, see for details
\cite{D1}, \cite{GWPL}.

If $\SS$ is a Whitney regular stratification of the reduced divisor $D_{red}$ and $x_0 \in D$ belongs to a (connected) stratum $S \in \SS$ with $d=\dim S>0$, then the local topology of the pair $(X,D)$ at the point $x_0$ is given by the product $(T,T\cap D)\times (S,0)$, where $(T,x_0)$ is a smooth transversal to $S$ at $x_0$, i.e. $\dim T=n-d$ and $T \pitchfork S$. In other words, the local topology of the pair $(X,D)$ is constant along the stratum $S$.

\noindent In terms of constructible sheaves, the topology of the pair $(X,D)$ is described by the vanishing cycle sheaf complex $\phi _f(\Q_X) \in D_c^b(\Q_D)$, see for instance \cite{D2}. If $i_T:T \to X$ denotes the inclusion of the transversal $T$ above, then there is an isomorphism
\begin{equation} \label{eq=vc1}
\phi _{f\circ i_T}(\Q_T)=i_T^{-1}\phi _f(\Q_X )
\end{equation}
in the derived category $ D_c^b(\Q_{T\cap D})$, see for instance \cite{Sch}, Lemma 4.3.4, p. 265. In fact this isomorphism holds for any $\SS$-constructible complex $\CC \in D_c^b(\Q_X)$, namely 
\begin{equation} \label{eq=vc2}
\phi _{f\circ i_T}(i_T^{-1}\CC   )=i_T^{-1}\phi _f(\CC ).
\end{equation}
This more general view-point allows us, in particular, to reduce the case of an arbitrary divisor $D$ to the case of a smooth hypersurface $H$, via the following basic construction. Set $X'=X \times \C$ and let $i_f: X \to X'$ be the graph embedding $x \mapsto (x,f(x))$. If $t$ denotes the coordinate on $\C$, then one has
\begin{equation} \label{eq=vc3}
\phi _{t}(i_{f*}\CC   )=i_{f*}\phi _f(\CC ).
\end{equation}
Indeed, $i_f$ is proper and $\phi _f$ commutes with proper direct images, see for instance \cite{D2}, p. 109.

\medskip

In this paper we prove two analogs of the above well-known properties.

\begin{thm} \label{thm=1} Let $D:f=0$ be a smooth divisor in $X$ and let $i_T: T \to X$ be the inclusion of a closed submanifold which is transversal to $D$. Assume that $M$ is a regular holonomic $\D_X$-module such that $T$ is non-characteristic for $M$ and for $M(*D)$. Let $V$ denote the Kashiwara-Malgrange filtration
of $M$ along the hypersurface $D$ and also the Kashiwara-Malgrange filtration of the restriction $i_T^*M$  along the hypersurface $T \cap D$. Then, for any $\alpha \in \C$, the following hold:

\medskip

\noindent (i)  $ i_T^*(V^{\alpha}M)=V^{\alpha}i_T^*M$;\,\,
(ii) $ i_T^*(Gr_V^{\alpha}M)=Gr_V^{\alpha}i_T^*M$.

\medskip\noindent In particular, 
$ i_T^*(\psi_fM)=\psi_{f|T}( i_T^* M)$ and $ i_T^*(\phi_fM)=\phi_{f|T}( i_T^* M)$.

\end{thm} 

For more details on the notions involved in this theorem we refer the reader to the next section. We mention here that the non-characteristic condition is a generalization of the transversality condition used in \eqref{eq=vc1} and  \eqref{eq=vc2} which implies that the analytic pull-back $i_T^*(M)$ coincides with the left derived pull-back $\bL i_T^*(M)$, see Proposition \ref{prop1}. Moreover, the V-filtration is the necessary ingredient to construct a theory of vanishing cycles for the regular holonomic $\D_X$-modules, see \cite{M},
 \cite{K2}, \cite{SM}, \cite{MM} for a synthesis. For other applications of Theorem \ref{thm=1} to the monodromy and Bernstein polynomials of families of hypersurface singularities, see  \cite{MT2}.

To state the second result, we return to the general case, i.e. $D$ is an effective divisor on $X$ defined by a holomorphic function $f$. The family of multiplier ideals $\{\I(\alpha D)\}_{\alpha \in \Q}$ associated to $D$ is a decreasing family of ideals in the structure sheaf $\OO_X$ which encodes a lot of the algebraic geometry of the pair $(X,D)$, see \cite{De}, \cite{K}, \cite{L}  for more on this beautiful and very active subject. Using the deep relation between multiplier ideals and V-filtrations established in \cite{BS} and Theorem \ref{thm=1} above, we obtain the following.

\begin{thm} \label{thm=2} Assume $T$ is transversal to any stratum of a Whitney regular stratification of the reduced divisor $D_{red}$, or more generally, $T$ is non-characteristic for the regular holonomic $\D_X$-module $\OO_X(*D)$. Then for $\alpha \in \Q$, we have a canonical isomorphism
$$i_T^*\I(\alpha D)=\I(\alpha (D \cap T)),$$
compatible with the inclusions $\I(\alpha D) \to \I(\alpha' D)$ and $\I(\alpha (D \cap T)) \to \I(\alpha' (D \cap T))$ for $\alpha > \alpha' $, where the isomorphism for $\alpha \le 0 $ is the natural isomorphism $i_T^* \OO_X =\OO_T$.
\end{thm}

Note that the multiplier ideals $\I(\alpha D)$ are essentially given by the filtration $\J(\alpha)$ on $\OO_X=\OO_X \otimes 1$ induced by the V-filtration on the $\D_{X'}$-module $B_f:=i_{f+}\OO_X=\OO_X \otimes \C[\partial _t]$ along the smooth hypersurface $H:t=0$, see  \cite{BS}. Here $i_{f+}$ denotes the direct image as a $\D$-module. Therefore Theorem  \ref{thm=2} is not a straightforward
consequence of Theorem \ref{thm=1} above. We can generalize Theorem \ref{thm=2} to the case of arbitrary subvarieties (see \ref{genth2}).
From Theorem  \ref{thm=2} we can deduce the following application to the spectrum in the sense of J.~Steenbrink \cite{St} (see (4.4-5)):

\begin{cor} \label{cor1}
Under the assumption of Theorem \ref{thm=2}, the spectrum $\Sp(f,x)$ of $f$ at $x\in T$ coincides with $\Sp(f|_T,x)$ up to a sign.
In particular, $\Sp(f|_T,x)$ is independent of $T$ in Theorem \ref{thm=2}.
\end{cor}

It has been known that the spectrum is
constant under a $\mu$-constant deformation of isolated hypersurface
singularities, see \cite{Va}. We have a weak generalization as follows.

\begin{cor} \label{cor11}
Let $S$ be a (connected) stratum of a Whitney regular stratification of $D_{red}$. Then $\Sp(f,x)$ and $\Sp(f|_T,x)$ are independent of $x\in S$, where $T$ is transversal to $S$ at $x$.
\end{cor}

Note that the $\mu$-constantness is equivalent to the Thom $a_f$-condition (see \cite{LS} and \ref{rk4} below) and the latter is weaker than the Whitney $(b)$ condition, see \cite{BMM1}. It is not clear whether Corollary \ref{cor11} holds assuming only the $a_f$-condition without the Whitney $(b)$ condition.

\section{Basic facts on $\D$-modules} \label{s1}

\subsection{Non-characteristic inverse images} \label{s11}

Let $X$  be a complex $n$-dimensional manifold and denote by $T^*X$ its cotangent bundle.
A point $(x,\xi) \in T^*X$ is just a pair formed by a point $x \in X$ and a linear form
$\xi:T_xX \to \C$, where $T_xX$ denotes the tangent space to $X$ at the point $x$.
When $Z \subset X$ is a locally closed analytic subset of $X$, we denote by $T^*_ZX$ the {\it conormal space} of $Z$ in $X$. This is by definition the closure in $T^*X|Z:=\pi ^{-1}(Z)$, with $\pi: T^*X  \to X$ the canonical projection, of the set of pairs
$(x,\xi) \in T^*X$ with $x \in Z$ a smooth point on $Z$ and $\xi|T_xZ=0$.
In particular, $T^*_XX$ is just the zero section of the cotangent bundle.

\begin{definition} \label{def1} Let $M$ be a coherent $\D_X$-module and let $CV(M) \subset T^*X$ be its characteristic variety. A submanifold $Z  \subset X$ is non characteristic for $M$ if
$$CV(M) \cap T^*_ZX \subset T^*_XX.$$

\end{definition}

The following basic example explains the relation to the transversality discussed briefly in the Introduction.

\begin{ex} \label{ex1} The Riemann-Hilbert correspondence, see \cite{M3},\cite{K4},  says that the DR-functor establishes an equivalence of categories
$$DR: D_{rh}^b(\D_X) \to  D_c^b(\C_X)$$
such that, for a regular holonomic $\D_X$-module $M$, the sheaf complex $\F=DR(M)$ is a perverse sheaf. Moreover, the characteristic variety $CV(M)$ coincides to the characteristic variety $CV(\F)$, which is defined topologically,  see for instance \cite{KS} or \cite{D2}, p. 111--113.
If $\SS$ is a Whitney regular stratification of $X$ such that $\F$ is $\SS$-constructible, then
$$CV(\F) \subset \cup_{S \in\SS}T_S^*X$$
 see for instance \cite{D2}, p. 119. It follows that a submanifold $T$ which is transversal to $\SS$, i.e. it is transversal to all the strata $S \in \SS$, it is automatically non-characteristic for $M$.
Note  that transversality of $T$ to a single stratum $S$ at a point $x \in S$ implies, via the (a)-regularity condition, transversality to all strata in a small neighborhood of $x$ in $X$.

To be even more specific, when $M=\OO_X(*D)$, then  $DR(\OO_X(*D))=Rj_*\C_U[n]$, where $U=X \setminus D$ and $j:U \to X$ is the inclusion, see \cite{Gro},  \cite{M4}. If $\SS$ is a Whitney regular stratification of $D$, then consider the Whitney regular stratification of $X$ given by   $\SS_0=\SS \cup \{U\} $. 
It follows that $Rj_*\C_U[n]$ is $\SS_0$-constructible and hence, if $T$ is transversal to $\SS$ as above, then $T$ is non-characteristic for $\OO_X(*D)$.

\end{ex}

The following result tells us that the property of being non-characteristic is preserved under small deformations.

\begin{lem} \label{lem1}
Let $M$ be  a  holonomic $\D_X$-module and $T \subset X$ a submanifold which is non-characteristic for $M$. Let $a \in T$ be any point in $T$ and $p:(X,a) \to (S,0)$ be the germ of a submersion onto a smooth germ $(S,0)$ such that the space germ induced by the special fiber $p^{-1}(0)$ coincides to $(T,a)$. Then there is an open neighborhood $U$ of $a$ in $X$ on which $p$ is defined and such that all the fibers $p^{-1}(s) \cap U$ for $s \in S$ are non-characteristic for the $\D_U$-module
$M|U$.

\end{lem}

\proof Since the question is local on $X$, we may suppose that $X=S \times T$ and $p$ is the first projection. On the other hand, it is known that the characteristic variety of a holonomic $\D_X$-module $M$ has the following local decomposition
\begin{equation} \label{eq=cv1}
CV(M)=\cup_{j=1,m}T^*_{Z_j}X
\end{equation}
where $Z_j$ are closed irreducible analytic subsets in $X$ and the conormal spaces $T^*_{Z_j}X$ are exactly the irreducible components of the characteristic variety $ CV(M)$. If a neighborhood $U$ as claimed above does not exist, then there is an index $j$ and a sequence of points $a_n \in X$ such that

\noindent (i) $a_n \to a$;

\noindent (ii) for each $n$, there is a point $(a_n,\xi _n) \in T^*_{Z_j}X$ such that $\xi_n \ne 0$ and
$\xi_n |0 \times T_{q(a_n)}T=0$, where $q:S \times T \to T$ is the second projection.

\noindent Moreover, we can norm the linear form $\xi_n$, e.g. by dividing by their norm and passing to a convergent subsequence, such that we may assume that

\noindent (iii) $\xi_n \to \xi \ne 0$.

\noindent Since $T^*_{Z_j}X$ is a closed subset, it follows that $(a,\xi)=\lim (a_n, \xi_n) \in 
T^*_{Z_j}X$. From (ii) we infer that $\xi |0 \times T_{q(a)}T=0$. Therefore
 $(a,\xi) \in T^*_TX\cap CV(M)$, in contradiction to the fact that $T$ is non-characteristic for $M$.

\endproof

We recall now the following result on non-characteristic inverse images, see for instance \cite{K1} or \cite{MT}, Prop. II.1.3, Thm. II.1.7.

\begin{prop} \label{prop1} Let $T$ be a submanifold in $X$ given by global equations $z_1=...=z_c=0$, where $c =\codim T$ such that $T$ is non-characteristic for a coherent $\D_X$-module $M$.
If $i_T:T \to X$ is the closed inclusion of $T$ into $X$, then the derived inverse image complex 
$ i_T^+M=  \bL i_T^*M$ is concentrated in degree zero and coincides to the  coherent $\D_T$-module 
$$i_T^*M=\frac{i^{-1}_TM}{z_1i^{-1}_TM+... +z_ci^{-1}_TM}.$$
Moreover, $z_1,...,z_c$ is a regular sequence in ${i^{-1}_TM}$ and $i_T^*M$ is holonomic when $M$ is holonomic.
More precisely, in this last situation, if $CV(M)=\cup_{j=1,m}T^*_{Z_j}X$, then 
$$CV(i_T^*M)=\cup_{j=1,m}T^*_{T \cap Z_j}T.$$
\end{prop}

The idea of the proof of this result is to determine the characteristic variety $CV(i_T^*M)$, see Thm. II.1.7 in  \cite{MT}. In the proof of Prop. II.1.3 in \cite{MT} it is shown that $\bL i_T^*M$ is represented by the Koszul complex of the sequence $z_1,...,z_c$ in $M$. Then this sequence is shown to be regular, by showing that any local section germ in $M$ is killed by some special type differential operators, similar to the ones we construct below in Lemma \ref{lem2}. See also
Corollary I.3.3 in \cite{MT}. It follows that the  Koszul complex has non-trivial homology only in degree zero, which yields the result.

\subsection{V-filtrations, b-polynomials and vanishing cycles} \label{s12}

Let $X$ be a complex analytic manifold and $H \subset X$ a smooth hypersurface. Let $\I$ be the ideal sheaf defining $H$.
We define the {\it decreasing Kashiwara-Malgrange V-filtration of $\D_X$ along $H$} by setting, for  $k \in \Z$,
\begin{equation} \label{eq=vf1}
V^k(\D_X)=\{P\in \D_X ~ | ~ P(\I^j) \subset \I^{j+k} \text { for all } j\in \Z\}
\end{equation}
where $\I^j=\OO_X$ for any $j<0$. It is easy to check that $V^k(\D_X) \cdot V^{\ell}(\D_X) \subset  V^{k+ \ell}(\D_X)$. In particular, $ V^{0}(\D_X)$ is a coherent sheaf of rings, see \cite{MM}, Prop. 2.1.5.

If $t=0$ is a local equation for $H$, then the differential operator $t\partial_t$ induces an element $E$ in $Gr^0_V\D_X=V^0(\D_X)/V^{1}(\D_X)$. This is called the Euler operator of $H$ and it is independent of the choice of the local equation  $t=0$. 

\medskip

A coherent $\D_X$-module $M$ is said to be {\it specializable along the hypersurface $H$} if it satisfies the following 
condition:

\medskip
\noindent ($\star$) For any point $x \in H$ and any germ $m\in M_x$, there is a non zero polynomial $b(s) \in \C[s]$ such that 
$b(E+1)m \in V^{1}(\D_X)_xm$.
\medskip

For more details, see \cite{MM}, Propositions II.1.9, II.2.2 and II.2.4.
 It is known that a holonomic $\D_X$-module is specializable along  any smooth hypersurface, see \cite{K2}. 
The {\it Bernstein polynomial} (or the {\it b-function}) $b_m$ of the germ $m\in M_x$ is the unitary polynomial of minimal degree satisfying the condition ($\star$). Note that $E+1=\partial _t \cdot t$. This shift by 1 is justified by the formula \eqref{eq=vf2} below.

\medskip

Let $<$ be a total order on $\C$ such that, for any $u,v \in \C$ we have $u<u+1$, $u<v$ if and only if $u+1<v+1$ and, finally,
there is some $m \in \N$ such that $v<u+m$. For instance, we can take the lexicografic order on $\C=\R^2$. Using this order, we define the {\it decreasing Kashiwara-Malgrange V-filtration along $H$} on the coherent $\D_X$-module $M$, assumed to be specializable along $H$, by
\begin{equation} \label{eq=vf2}
V^{\alpha}M_x = \{m \in M_x ~|~ \text{ all the roots of the $b$-function}~~ b_m ~~\text{are}~ \geq \alpha \}.
\end{equation}
See Kashiwara \cite{K2}, Malgrange \cite{M}, and also \cite{Sab}. The filtration $V^{\alpha}M$ is indexed by a finite union of lattices $\beta +\Z$ in $\C$, hence it is a discrete, decreasing and exhaustive filtration on $M$.
In most cases coming from geometry, e.g. when $M$ is obtained by applying some natural functors to the 
$\D_X$-module $\OO_X$, then all $\beta \in \Q$, and hence $V$ is actually indexed by $\Q$. This is the case in particular for the module $M=B_f$.

\noindent It turns out that the V-filtration can be defined by the following list of characteristic properties,  see \cite{Sab}, \cite{S3} for more details. Define $V^{> \alpha}M =\cup_{\beta >\alpha}  V^{\beta}M$ and $Gr_V^{ \alpha}M=V^{\alpha}M/V^{>\alpha}M$.

\begin{prop} \label{prop3}  Let $M$ be a coherent $\D_X$-module, specializable along $H$. 
The Kashiwara-Malgrange V-filtration is the unique discrete, decreasing and exhaustive filtration on $M$ satisfying the following conditions: 

\begin{enumerate}

\item  $V^k(\D_X) \cdot V^{\alpha}M \subset  V^{k+ \alpha}M$, for any $k \in \Z$ and $\alpha \in \C$;

\item $V^{\alpha}M$ is a coherent $V^0(\D_X)$-module, for any  $\alpha \in \C$;

\item $t \cdot V^{\alpha}M=V^{\alpha +1}M$, for  $\alpha >> 0$;

\item the action of $\partial _t \cdot t - \alpha$ on $Gr_V^{ \alpha}M$ is locally nilpotent.

\end{enumerate}

\end{prop}

\begin{rk} \label{rk1} \rm
Many authors prefer to work with the {\it increasing} V-filtration defined by $V_{\alpha}M=V^{-\alpha}M$, see for instance
\cite{MM}. Our choice here is justified by the fact that, as we mentioned in the Introduction, in an important case the
V-filtration is essentially the filtration induced by the family of multiplier ideals, which is by definition a decreasing filtration.
\end{rk}

\begin{definition} \label{def3} Let $M$ be a regular holonomic $\D_X$-module. Let
$t=0$ be a global equation for the hypersurface $H$. Then the nearby and the vanishing cycle functors $\psi_t, \phi_t: Mod_{rh}(\D_X) \to
 Mod_{rh}(\D_H)$ are defined as follows:
$$\psi_tM=\oplus_{0< \alpha \leq 1}Gr_V^{ \alpha}M$$
and
$$\phi_tM=\oplus_{0\leq  \alpha <1}Gr_V^{ \alpha}M.$$
\end{definition}
These functors are related to the topological perverse nearby and  vanishing cycle functors  ${}^p\psi_t,{}^p \phi_t: Perv(X,\C) \to  Perv(H,\C)$, see for instance
\cite{D2}, p. 139, by the following relations $DR \circ \psi_t={}^p\psi_t \circ DR$ and  $DR \circ \phi_t={}^p\phi_t \circ DR$.
We recall also the exact triangle in $D^b(\D_H)$
\begin{equation} \label{eq=vc4}
 \psi_tM \to \phi_tM \to  i_H^+M  \to 
\end{equation}
where $i_H:H \to X$ is the inclusion and the morphism $\can: \psi_tM \to \phi_tM $ is induced by $\partial_t$.

In the sequel we will regard not only $\psi_tM$, $\phi_tM$ but also any $Gr_V^{ \alpha}M$ as $\D_T$-modules.
Let  $\Psi_tM= i_{H+} \psi_tM$ and  $\Phi_tM= i_{H+} \phi_tM$ be the corresponding  $\D_X$-modules supported by the hypersurface $H$.

\subsection{Relatively specializable $\D_X$-modules and relative V-filtrations} \label{s13}

In the situation of the previous subsection, suppose that we have in addition a submersion $p:X \to S$ such that the composition
$p_H=p\circ i_H$ is still a submersion. We denote by $\D_{X/S} \subset \D_X$ the sheaf of relative differential operators, which is by definition the sheaf of subrings in $ \D_X$ spanned by $\OO_X$ and by the derivations coming from vector fields on $X$, tangent to the fibers of $p$. We define the relative V-filtration on  $\D_{X/S}$ simply by setting $V^k(\D_{X/S})=V^k(\D_X)\cap \D_{X/S}$.

\begin{definition}
 A coherent $\D_{X/S}$-module $M$ is relatively specializable along $H$ if there exists a decreasing, exhaustive
filtration $U^k(M)_{k\in\Z}$ of $M$ by coherent $V^0(\D_{X/S})$-submodules such that the following conditions
are satisfied:

\begin{enumerate}
\item $V^k(\D_{X/S}) \cdot U^\ell(M) \subset  U^{k+\ell}(M)$, for any $k,\ell \in \Z$;

\item Locally on $X$, there is an $\ell \in \N$ such that,  for any integer $k \in \N$, we have
 $V^k(\D_{X/S}) \cdot U^\ell(M)=  U^{k+\ell}(M)$ and $V^{-k}(\D_{X/S}) \cdot U^{-\ell}(M)=  U^{-k-\ell}(M)$;

\item Locally on $X$, there exists a nonzero polynomial $b(s) \in \C[s]$
such that $b(E-k)U^k(M) \subset U^{k+1}(M)$, for all $k\in \Z$.
\end{enumerate}
\end{definition}

For more details, see \cite{MT2}. In fact, one can see that if $M$ is a  coherent 
$\D_{X}$-module, specializable along $H$, and such that the $V^0(\D_X)$-modules
 $V^\alpha(M)$ are coherent 
over $V^0(\D_{X/S})$, then $M$ is relatively specializable along $H$
as a $\D_{X/S}$-module.
Moreover,
the converse implication is easy to check. In the following, we will use this last 
characterization.

\subsection{Characteristic varieties} \label{s14}

Let $X$ be a complex $n$-dimensional manifold, $Z \subset X$ an irreducible analytic subset and $f:X \to \C$ a holomorphic function such that the restriction $f|Z$ is not a constant function, i.e. $Z$ is not contained is a fiber of $f$.
\begin{definition} \label{def4}
The relative conormal space of the restriction  $f|Z$ is the closed analytic subset in $T^*X$ given by
$$T^*_{f|Z}X=\overline {\{(x,\xi+\lambda df(x))~|~ (x,\xi)\in T_Z^*X, ~ \lambda \in \C\}}$$
where the closure is taken in  $T^*X$. The associated Lagrangian variety  of the restriction  $f|Z$ is the subset in
$T^*X|f^{-1}(0)$ given by
$$W_0({f|Z})=T^*_{f|Z}X\cap \pi^{-1}(f^{-1}(0))$$
where $\pi:T^*X \to X$ is the canonical projection, see  \cite{BMM1}, \cite{Gin}.
\end{definition}
It follows from \cite{K0},  \cite{Gin}, Prop. 2.14.1 that $T^*_{f|Z}X$ is an $(n+1)$-dimensional subvariety in $T^*X$ and
$W_0({f|Z})$ is a closed, conic, Lagrangian subvariety in $T^*X$; in particular, $\dim W_0({f|Z})=n$. 

Let $\ZZ=(Z_a)_{a \in A}$ be an analytic stratification of the analytic set $Z$, i.e. $Z=\cup_{a \in A}Z_a$ is a partition of $Z$ into smooth semianalytic subsets $Z_a$ satisfying the frontier condition: if $Z_a \cap {\overline Z}_b\ne \emptyset$, then  $Z_a \subset {\overline Z}_b$. Similarly, let $\T=(T_b)_{b \in B}$ be an analytic stratification of $\C$, such that the pair $(\ZZ,\T)$ give a stratification for $f:Z \to \C$, i.e. for any $a \in A$ there is a $b\in B$ such that $f(Z_a)\subset T_b$ and the induced mapping $f:Z_a \to T_b$ is a submersion.

\begin{definition} \label{defaf}
We say that the stratification  $(\ZZ,\T)$ of $f:Z \to \C$ satisfies the Thom $a_f$-condition if for any pair of strata  $Z_a \subset {\overline Z}_b$ one has
$$T^*_{f|{\overline Z}_b}X\cap\pi^{-1}(Z_a) \subset T^*_{f|Z_a}X$$
where we set $ T^*_{f|S}X=T^*_SX$ whenever $f|S$ is a constant function.
\end{definition}

In particular, if $\ZZ$ is a stratification of the pair $(Z,Z_0)$ where $Z_0= f^{-1}(0) \cap Z$, it follows that the induced stratification $\ZZ_0$ on $Z_0$ satisfies the Whitney $(a)$-condition. In such a case 
\begin{equation} \label{eqaf}
 W_0(f|Z) =\cup _S( T^*_{f|Z}X\cap\pi^{-1}(S)) \subset \cup _S T^*_SX
\end{equation}
where $S$ runs through the strata of $\ZZ_0$. This last union is a closed subset, since this property is equivalent to the  Whitney $(a)$-condition for the induced stratification  $\ZZ_0$.
See also \cite{HMS} for other results about the relative conormal space
and the relation with Thom's ($a_f$)-condition. 

We recall the geometric formulation of Ginsburg results in  \cite{Gin}, Prop. 2.14.4 and Thm. 3.3 and 5.5 given by Brian\c con, Maisonobe and Merle in  \cite{BMM1},
Thm 3.4.2.

\begin{prop} \label{prop5}
Let $M$ be a regular holonomic $\D_X$-module, $M_0\subset M$ a coherent $\OO_X$-submodule spanning $M$ over $\D_X$,
and $f:X \to \C$  a holomorphic function. Assume that the characteristic variety of $M$ is given by
$CV(M)=\cup_{j}T^*_{Z_j}X$ as in \eqref{eq=cv1}. Then the following equalities hold, the first two ones 
in a neighborhood of the zero set $f=0$.
\begin{enumerate}

\item $CV(\D_X[s]M_0f^s)=\bigcup_{f|Z_j \ne 0}T^*_{f|Z_j}X$;

\item $CV(M[1/f])=\bigcup_{f|Z_j \ne 0}T^*_{Z_j}X \cup \bigcup_{f|Z_j \ne 0}W_0(f|Z_j)$;

\item If $H:f=0$ is a smooth hypersurface, then by setting $t=f$, we have 
$CV(\Psi_t M)= \bigcup_{f|Z_j \ne 0}W_0(f|Z_j)$
in $T^*X$.

\end{enumerate}

\end{prop}

Here the  $\D_X$-module $\D_X[s]M_0f^s$ is defined as follows. Consider the  $\OO_X$-submodule
$M_0f^s$ in the  $\D_X$-module $M[1/f,s]f^s=  M \otimes _{\OO_X}\OO_X[1/f,s]f^s$, where the action of
 $\D_X$ on  $M[1/f,s]f^s$ is the usual one, e.g.
$$\frac{\partial}{\partial x_k}\cdot f^s=s \cdot \frac{\partial f}{\partial x_k}\cdot 1/f \cdot f^s.$$
Then  $\D_X[s]M_0f^s$ is the  $\D_X$-module spanned by  $M_0f^s$ in  $M[1/f,s]f^s$, which is known to be coherent
(for the case $M_0=\OO_X \cdot m$ see \cite{MT}, Lemme III.1.2).

\medskip

Concerning the last property in Prop. \ref{prop5}, note that $\bigcup_{f|Z_j \ne 0}W_0(f|Z_j)=\cup _kT^*_{Y_k}X$, where $(Y_k)_k$ is a locally finite family of irreducible subvarieties in $H$. It follows that, when we regard $\psi_t M$, one has the decomposition
\begin{equation} \label{eq=cv3}
 CV(\psi_t M)= \cup _kT^*_{Y_k}H.
\end{equation}

\begin{cor} \label{cor2}
Let  $X$ be a complex  manifold, $H \subset X$ a smooth hypersurface given by $t=0$ and $p:X \to S$ a submersion such that
the restriction $p_H=p|H$ is still a submersion. Let  $M$ be a regular holonomic $\D_X$-module such that the fibers of $p$ are non-characteristic for $M$ and $M(*H)$. Then the fibers of the restriction $p_H$ are non-characteristic for
the nearby cycle module $\psi_tM$, the vanishing cycle module $\phi_tM$ and any  $\D_H$-module $Gr_V^{\alpha}M$, for $
\alpha \in \C$.

\end{cor}

\proof

Let $F_s=p^{-1}(s)$ be a fiber of $p$ and note that $F_s$ is transversal to $H$. Then, using a similar formula to
\eqref{eq=cv3}, we infer that,
for a $\D_H$-module $N$, 
$F_s \cap H$ is non-characteristic for $N$ if and only if  $F_s $ is non-characteristic for $i_{H+}N$. 
 Using Prop. \ref{prop5}, we get
$$CV(\Psi_t M) \subset  CV(M(*H)).$$
This implies that  $F_s $ is non-characteristic for $\Psi_t M$ and hence $F_s \cap H$ is non-characteristic for
 $\psi_t M$.
On the other hand, the exact triangle
\begin{equation} \label{eq=tr1}
 R\Gamma_H(M)_{alg} \to M \to M(*H) \to
\end{equation}
implies $CV( R\Gamma_H(M)_{alg}       )\subset CV(M)~ \cup ~ CV(M(*H))$, and hence 
 $F_s $ is non-characteristic for $ R\Gamma_H(M)_{alg}$. Since $H$ is smooth, it follows that
$$R\Gamma_H(M)_{alg}=i_{H!}i_H^!(M)$$
see \cite{DMSS}, 1.5.4-1.5.5. Using now the fact that a regular holonomic $\D_X$-module $M$ and its dual $DM$
have the same characteristic variety and that $D(i_{H!}i_H^!(M))=i_{H+} i_H^+M$, it follows that
 $F_s $ is non-characteristic for $i_{H+} i_H^+M$. Using now the exact triangle \eqref{eq=vc4},
 it follows that $F_s \cap H$ is non-characteristic for $\phi_tM$.

Since $CV(M_1\oplus M_2) \supset CV(M_1)$, it follows from the above that  
$F_s \cap H$ is non-characteristic for  $Gr_V^{\alpha}M$, for $
\alpha \in \C$, $0 \leq \alpha \leq 1$. The proof is completed using the fact that
$$t: Gr_V^{\alpha}M \to Gr_V^{\alpha+1}M$$
is an isomorphism of $\D_H$-modules, for any $\alpha \ne 0$.

\endproof

\begin{rk} \label{rk2} \rm
The proof above shows that, with the above notation, for a submanifold $T$ in $X$, the condition  

\medskip

\noindent (C1) $T$ is non-characteristic for  $M$ and $M(*H)$

\medskip

is equivalent to the condition

\medskip

\noindent (C2) $T$ is non-characteristic for  $M$ and $i_{H+} i_H^+M$.

\end{rk}

Lemma \ref{lem1} and Corollary  \ref{cor2} yield the following

\begin{cor} \label{cor3}
Let  $X$ be a complex  manifold, $H \subset X$ a smooth hypersurface given by $t=0$ and $T \subset X$ a submanifold which is transverse to $H$. Let  $M$ be a regular holonomic $\D_X$-module such that  $T$ is non-characteristic for $M$ and $M(*H)$. Then the submanifold $T \cap H$ is non-characteristic for
the nearby cycle module $\psi_tM$, the vanishing cycle module $\phi_tM$ and any  $\D_H$-module $Gr_V^{\alpha}M$, for $
\alpha \in \C$.

\end{cor}

\section{Proof of Theorem \ref{thm=1}} \label{s2}

By writing $T$ as the intersection of $d$ smooth hypersurfaces, each of them transversal to $D$, we may restrict our
attention to the case $d=\dim S= \codim T =1$.
We start with the following general result, giving a sufficient condition for a regular holonomic module to be relatively specializable. Recall the setting and the notation from section \ref{s13}.

\begin{prop} \label{prop6}

Let $M$ be a regular holonomic $\D_X$-module such that the fibers of $p:X \to S$ are non-characteristic for $M$ and for $M(*H)$. Then $M$ is relatively specializable along the smooth hypersurface $H \subset X$.

\end{prop}

\proof

Let $a$ be a point in $H$ and choose a local system of coordinates 
$$(x_1, \ldots ,x_{n-2}, y,t)$$
 at $a$ in $X$ such that $h(x_1, \ldots ,x_{n-2}, y,t)  = t=0$ is a local equation for $H$ and
$p(x_1, \ldots ,x_{n-2}, y,t) = y$. To proceed,
we need the following.

\begin{lem}\label{lem2} 
With the above assumptions, for any local section germ $m \in M_a$ or $m \in M(*H)_a$, there is a differential operator
$\widetilde{P}\in V^0({\D}_{X})$ 
killing $m$ and having the following special form
$$\widetilde{P} = \left( \frac{\partial}{\partial y} \right) ^{k} + A_{1}  \left( \frac{\partial}{\partial y}\right) ^{k -1} + \cdots + A_{k}\ ,$$
where $A_{i}$ are  operators of degree  $\leq i$,
independent of the derivatives ${\partial}/{\partial t}$, ${\partial}/{\partial y}$.
 
\end{lem}

\proof

Let ${M}_0 \subset {M}(*H)$ be a coherent $\OO_X$-module spanning $M(*H)$ over $\D_X$.
It follows then from Proposition \ref{prop5}, (i) that the fibers of $p$ are  non-characteristic for the
coherent ${\D}_{X}$-module ${\D}_{X}[s] { M}_0  h^s$. In particular,
for any germ $m \in  M(*H)$, there is a differential operator
 $P$ of total degree $k$, killing $mt^s$,  and of the form
$$ P = \left( \frac{\partial}{\partial y} \right) ^{k} + A_{1}  \left( \frac{\partial}{\partial y}
\right)^{k-1} + \cdots + A_{k}\  .$$
Here the  differential operators  $A_{i}$ are
of degree $ \leq i$,
independent  of the  derivation ${\partial}/{\partial y} $ 
(compare to the proof of Proposition II.1.3 in \cite{MT}). 

Recall also the simple fact that an equality
$$\sum_{i=0}^k\left(\frac{\partial}{\partial t}\right)^i m_i t^s= 0$$
with $m_0,\ldots,m_k \in{ M}(*H)_a$ implies 
$m_0=...=m_k=0$, see for instance \cite{MM} Lemme 2.4.1.
Therefore, there is a  differential operator  $\widetilde{P}$, 
 having the same properties as the operator
$P$, and being in addition independent of ${\partial}/{\partial t}$. 

To end the proof of this Lemma, we have still to consider the case of a germ
 $m\in{ M}$ killed by a power of $t$. In other words, we may suppose now that ${ M}$ is supported by
 the hypersurface $H$. In particular, ${ M}$ is the image under $i_{H+}$ of a ${\D}_H$-module ${N}$ 
such that the fibers of  $p_T=p \circ i_T$ are  non-characteristic for    
${N}$. But then the same argument as above applied to $N$ shows that the germs of sections in ${N}$, and hence in  ${ M}$, are killed by  a differential operator
similar to the operator  $\widetilde{P}$, 
but which in this case are independent of the variable $t$ and of the derivative ${\partial}/{\partial t}$. 
\endproof

Let us return now to the proof of Proposition \ref{prop6}. Let  $m\in{M}_a$ be the germ of a local section. Using the  corresponding differential operator     $\widetilde{P}
\in V^0({\D}_X)$ obtained as in the above Lemma, 
we get
$$ V^\ell({\D}_X )_a m = \sum_{\kappa < k} V^\ell({\D}_{{X} \mid {S}} )_a
({\partial}/{\partial y})^{\kappa}  m $$
for all $\ell \in {\Z}$. It follows that the $V^0(\D_X)$-modules $V^\alpha(M)$ are
$V^0(\D_{X/S})$-coherents; in particular, the specializable module $M$ is
relatively specializable along $H$. 
\endproof

Finally we can prove Theorem \ref{thm=1}. Using Lemma \ref{lem1} we can place ourselves in the relative
case of a fibration $p:X \to S$ such that $T=p^{-1}(0)$ is the special fiber of $p$ and $D=H$ is the smooth hypersurface.
If we regard $ V^{\alpha}M$ as a $V^0\D_X$-module, then one has 
\begin{equation} \label{eq=Li}
\bL i_T^*( V^{\alpha}M)= i_T^*( V^{\alpha}M).
\end{equation}
This follows from a  $V^0\D_X$-version of  Proposition \ref{prop1}, using the differential operator constructed in Lemma
 \ref{lem2} to show that the action of $y$ is injective on $M$. Note that there is a natural morphism
\begin{equation} \label{eq=n}
\nu: i_T^*( V^{\alpha}M) \to i_T^*(M) 
\end{equation}
induced by the inclusion $ V^{\alpha}M \to M$. This morphism $\nu$ is actually injective. To see this,
apply the functor $\bL i_T^*$ to the exact sequence of $V^0\D_X$-modules
\begin{equation} \label{eq=es1}
0 \to  V^{\alpha}M \to M \to \frac {M}{ V^{\alpha}M} \to 0
\end{equation}
and note that the action of $y$ on the last term is injective by the same argument as above; in particular, we have
$$\bL i_T^*\frac {M}{ V^{\alpha}M}= i_T^*\frac {M}{ V^{\alpha}M}.$$
To prove that the image of the morphism $\nu$ is exactly $ V^{\alpha} i_T^* M$, we proceed as follows.
First, we get as above
$$\bL i_T^* V^0\D_{X|S}= i_T^* V^0\D_{X|S}=V^0(\D_T).$$
Then, applying the functor $ i_T^*=\otimes _{ V^0\D_{X|S}} V^0(\D_T)  $ to a presentation of $ V^{\alpha}M$ as a coherent $ V^0\D_{X|S}$-module, which exists via Proposition  \ref{prop6}, we get that
$\nu( i_T^*( V^{\alpha}M))$ is a coherent $V^0(\D_T)$-module. Since the filtration on $i_T^* M$
given by $(\nu  ( i_T^*( V^{\alpha}M)))_{\alpha}$ clearly satisfies all the other conditions in 
Proposition  \ref{prop3}, it follows that 
$$\nu  ( i_T^*( V^{\alpha}M))= V^{\alpha} i_T^* M$$
which ends the proof of the first claim in  Theorem \ref{thm=1}.

\medskip
The second claim in  Theorem \ref{thm=1} follows by applying the functor $\bL i_T^*= i_T^*$ to the exact sequence
of  $V^0(\D_X)$-modules
$$0 \to  V^{ >\alpha}M \to  V^{\alpha}M \to Gr_V^{\alpha}M \to 0.$$

\section{Proof of Theorem  \ref{thm=2} } \label{s3}

It is clearly enough to treat the case when $T$ is a hypersurface. Since the question is local, we will assume that $X=\C^n$, $x_0=0$ is the origin. Let us choose a system of coordinates $(y,z)$ at the origin, with
$y=(y_1,...,y_{n-1})$ such that $T:z=0$. Then the induced effective divisor $T \cap D$ is given by an equation $f_T(y)=f(y,0)=0$. Consider the following diagram of closed embeddings, where $T'=T\times \C$.

$$\xymatrix{
T \ar[r]^{i_{f_T}} \ar[d]^{i_T} & T' \ar[d]^ {i_{T'}} \\
X   \ar[r]^{i_f } & X'.
}$$
Let $t$ be a coordinate on $\C$ and set as in the Introduction 
$B_f=(i_f)_+\OO_X = \OO_X \otimes \C[\partial _t] $, where $\partial _t= \partial / \partial t$ and 
$\J(\alpha)=V^{\alpha}B_f \cap (\OO_X \otimes 1).$ Similarly, we set 
$B_{f_T}=(i_ {f_T} )_+\OO_T$  and $\J_T(\alpha)=V^{\alpha}B_{f_T} \cap (\OO_T \otimes 1).$
Using Theorem 0.1 in \cite{BS}, in order to prove Theorem \ref{thm=2}, it is enough to show
that
\begin{equation} \label{eq=mi1}
i_T^*\J(\alpha)=\J_T(\alpha)
\end{equation}
for all $\alpha \in \Q$.

By hypothesis $T'$ is non-characteristic for the $\D_{X'}$-modules $B_f$ and $B_f[\frac{1}{t}]$ (using the coordinate $t':=t-f$).
Proposition \ref{prop1} then yields
\begin{equation} \label{eq=mi2}
\bL i_{T'}^*B_f= i_{T'}^*B_f=B_{f_T}
\end{equation}
and we infer from Theorem \ref{thm=1} that
\begin{equation} \label{eq=mi3}
 i_{T'}^*V^{\alpha}  B_f= V^{\alpha} B_{f_T}
\end{equation}
for all $\alpha \in \Q$.
The module $B_f$ is endowed with a natural increasing Hodge filtration, given up-to a shift by
\begin{equation} \label{eq=hodge1}
 F_p B_f=\oplus _{0\leq j \leq p} \OO_X \otimes \partial _t ^j
\end{equation}
for all $p \in \Z$. In particular $ \OO_X \otimes 1=F_0B_f$.

Using the Hodge filtration $F$, we may write
$$\J(\alpha)=V^{\alpha}  B_f \cap F_0B_f$$
and similarly
$$\J_T(\alpha)=V^{\alpha}  B_ {f_T}\cap F_0B_{f_T}  .$$
Since obviously $i_{T'}^*(\OO_X \otimes \partial _t ^j)=\OO_T \otimes \partial _t ^j$, it follows that
\begin{equation} \label{eq=hodge2}
 i_{T'}^*(F_pB_f) \simeq F_pB_{f_T}
\end{equation}
for all $p \in \N$. Hence the relation \eqref{eq=mi1} above is equivalent to
\begin{equation} \label{eq=hodge3}
 i_{T'}^*( V^{\alpha}  B_f \cap      F_0B_f)  \simeq    i_{T'}^*( V^{\alpha}  B_f) \cap  i_{T'}^*(  F_0B_f).
\end{equation}
Therefore, we complete the proof of  Theorem  \ref{thm=2} if we prove the following

\begin{lem} \label{lem3}
The isomorphism
$$i_{T'}^*( V^{\alpha}  B_f \cap      F_pB_f)  \simeq  i_{T'}^*( V^{\alpha}  B_f) \cap  i_{T'}^*(  F_pB_f)$$
holds for  all $p \in \Z$ and  all $\alpha \in \Q$.
\end{lem}

Set $M=B_f$, $F_p V^{\alpha}M=F_pM \cap  V^{\alpha}M$, $V^{\alpha}(M/F_pM)= V^{\alpha}M/F_p V^{\alpha}M$
and $F_p(M/V^{\alpha}M)=F_pM/F_p V^{\alpha}M$. Then the following diagram, where all the morphisms are the
canonical monomorphisms or canonical epimorphisms, has obviously exact rows and exact columns.
\begin{equation} \label{eq=diag1}
\xymatrix{
        & 0 \ar[d]                      &0  \ar[d]                  &0  \ar[d]                        & \\
0\ar[r] & F_p V^{\alpha}M  \ar[r] \ar[d]& V^{\alpha}M  \ar[r] \ar[d]& V^{\alpha}(M/F_pM) \ar[r] \ar[d]& 0\\
0\ar[r] &  F_pM  \ar[r] \ar[d]          &           M  \ar[r] \ar[d]&            M/F_pM \ar[r] \ar[d]& 0\\
0\ar[r] &  F_p(M/V^{\alpha}M)\ar[r] \ar[d]& M/V^{\alpha}M\ar[r] \ar[d]& M/(F_pM+V^{\alpha}M)\ar[r] \ar[d]& 0\\
        & 0                               &0                          &0                                 &
}
\end{equation}

Lemma \ref{lem3} is clearly equivalent to the following

\begin{lem} \label{lem4}
The diagram 
$$\xymatrix{
        & 0 \ar[d]                      &0  \ar[d]                  &0  \ar[d]                        & \\
0\ar[r] & i_{T'}^*(F_p V^{\alpha}M)  \ar[r] \ar[d]& i_{T'}^*(V^{\alpha}M)  \ar[r] \ar[d]& i_{T'}^*(V^{\alpha}(M/F_pM)) \ar[r] \ar[d]& 0\\
0\ar[r] & i_{T'}^*( F_pM ) \ar[r] \ar[d]          &  i_{T'}^*(         M)  \ar[r] \ar[d]&  i_{T'}^*(          M/F_pM) \ar[r] \ar[d]& 0\\
0\ar[r] & i_{T'}^*( F_p(M/V^{\alpha}M))\ar[r] \ar[d]& i_{T'}^*(M/V^{\alpha}M) \ar[r] \ar[d]& i_{T'}^* (M/(F_pM+V^{\alpha}M)) \ar[r] \ar[d]& 0\\
        & 0                               &0                          &0                                 &
}$$
obtained from the above diagram by applying the functor $i_{T'}^*$ has still
exact rows and exact columns.

\end{lem}

Now, if we have an exact sequence of $R$-modules
$$0 \to N' \to N \to N'' \to 0$$
and if the action of $z \in R$ is injective on  $N''$, then it follows, e.g. by the use of Tor exact sequence, that the induced sequence
$$0 \to N'/zN' \to N/zN \to N''/zN'' \to 0$$
is still exact. It follows that all we have to prove is in fact the following.

\begin{lem} \label{lem5}
The action of $z \in \OO_{X'}$ is injective on any of the five $\OO_{X'}$ -modules
occurring in the third column and the third row of  the diagram \eqref{eq=diag1}.
\end{lem}

\proof

We have already proved that the action of $z$ is injective on $M/V^\alpha M$
(see the proof of Theorem \ref{thm=1}). On the other hand, the
assumption for $M/F_pM$ is clear since $M/F_pM$ is $\OO_{X'}$-free. Thus
we just have to check the assertion for $M/(F_pM+V^\alpha M)$.  
In fact we prove a slightly stronger claim, i.e. we show that the action of $z$ is injective on
any of the quotients
\begin{equation} \label{eq=claim1}
 F_p/F_q( V^{\alpha}/ V^{\beta})M=\frac{F_p V^{\alpha}M}{F_p V^{\beta}M+F_q V^{\alpha}M}
\end{equation}
for any $p>q$ in $\Z$ and any $\alpha < \beta$ in $\Q$.

\medskip

Here we can use the following well-known fact:

\medskip
Let $N_1\subsetneq N_2\subsetneq \cdots \subsetneq N_p$ be
an increasing sequence of $R$-modules. Assume that the action of $z\in R$
is injective on $N_i/N_{i-1}$ for $2\leq i\leq p$. Then the action of $z$
is injective on $N_i/N_j$ for $1\leq j<i\leq p$.

\medskip
Therefore it is enough to show that the action of $z$ is injective on 
$M_p^{\alpha}=Gr_p^FGr_V^{\alpha}M$. 
Consider the graded $Gr^F \D_{X}$-module
$$M^{\alpha}=\oplus _p M_p^{\alpha}$$
and note that
\begin{equation} \label{eq=supp1}
\supp M^{\alpha}=CV(Gr_V^{\alpha}M).
\end{equation}
Suppose the action of $z$ on $M^{\alpha}$ is not injective.
Then some irreducible component of $\supp M^{\alpha}$ is contained in $\{z=0\}$ since $M^{\alpha}$ is a Cohen-Macaulay $Gr^F \D_{X}$-module (see \cite{S1}, Lemme 5.1.13). Indeed, if $zm = 0 $ for some nonzero $m\in M^{\alpha}$, then the support of the submodule generated by $m$ is contained in $\{z=0\}$ and its dimension coincides with that of $M^{\alpha}$ by the Cohen-Macaulay property. On the other hand, $T$ is non-characteristic for $Gr_V^{\alpha}M$ by Prop. \ref{prop5} because $T'$ is non-characteristic for $\B_f[\frac{1}{t}]$. This is clearly a contradiction.
Thus we get the injectivity of the action of $z$.
This completes the proof of Theorem \ref{thm=2}.

\subsection{Spectrum} \label{spectrum}
For a holomorphic function $f$ on a complex manifold $X$ and
$x\in f^{-1}(0)$, the spectrum $\Sp(f,x)=\sum_{\alpha\in\Q}
m_{\alpha}t^{\alpha}$ is defined by$$\aligned&m_{\alpha} = \sum_{j} (-1)^{j-n+1} \dim Gr_{F}^{p}\tH^{j}(F_x,\C)_{\lambda}\\&\quad\quad\text{with}\,\, p = [n - \alpha],\,\,\lambda = \exp(-2\pi i\alpha),\endaligned$$where $n=\dim X$, $F_x$ denotes the Milnor fiber of $f$ around $x$,
$\tH^{j}(F_x),\C)_{\lambda}$ is the $\lambda$-eigenspace of the
reduced cohomology for the semi-simple part of the Milnor monodromy,
and $F^p$ denotes the Hodge filtration, see \cite{St} (and also \cite{S3}).

\subsection{Proof of Corollary \ref{cor1}} \label{cor1proof}
By the arguments in the proof of Theorem \ref{thm=2}, the three filtrations $F,V_t,V_z$ on $\B_f$ are compatible in the sense of \cite{S1}, 1.1.13. (Indeed, $V_z$ is the $z$-adic filtration in this case and applying the multiplication by $z^j$ to the diagram (\ref{eq=diag1}) in the proof of Lemma \ref{lem3} we get a cubic diagram of short exact sequences which is equivalent to the compatibility of the three filtrations. Here $V_t$ denotes the $V$-filtration along $t=0$ and similarly for $V_z$.)
This compatibility implies $i_T^*(\varphi_f(\OO_X,F))=\varphi_{f|T}(i_T^*(\OO_X,F))$, because $i_T^*=\psi_z$ and $\varphi_z=0$ on $(\varphi_f(\OO_X,F)$ in this case.

Let $i_x:\{x\}\to X$ denote the inclusion morphism and $i_x^*$ be the pull-back in the category of filtered $\D$-modules underlying mixed Hodge modules.
This is defined by iterating the mapping cone of
$\psi_{x_i,1}\to\varphi_{x_i,1}$ where the $x_i$ are local coordinates.
Then the Hodge filtration on the Milnor cohomology is given by $H^ji_x^*(\varphi_f(\OO_X,F))$, see also \cite{BM}.
So the assertion follows.

\subsection{Proof of Corollary \ref{cor11}} \label{cor11proof}
Let $i_S:S\to X$ denotes the inclusion. By Proposition \ref{prop5} the characteristic variety of the nearby cycle sheaf is contained in the union of the conormal bundles of the strata of the stratification satisfying the Thom $a_f$-condition (see \ref{defaf}). Therefore the $\cH^ji_S^*\psi_f\C_X$ for $j\in\Z$ are local systems by \cite{KS}, Proposition 8.4.1 (note that a $\mu$-stratification in loc.~cit. is Whitney regular, see Trotman \cite{Tr}). Since the pull-back by $i_S$ is compatible with the pull-back in the derived category of mixed Hodge modules, these local systems naturally underly variations of mixed Hodge structures. (Note that a mixed Hodge module is a variation of mixed Hodge structure if its underlying perverse sheaf is a local system up to a shift of complex). So the Hodge numbers are constant, and we get the assertion for $\Sp(f,x)$. Then the assertion for $\Sp(f|_T,x)$ follows from Corollary \ref{cor1}.

\begin{rk} \label{rk3} \rm
It is not clear whether Corollary \ref{cor11} holds assuming only the $a_f$-condition. Note that for $K\in D_c^b(\C_X)$ and a stratification satisfying only the Whitney $(a)$ condition and such that the characteristic variety $CV(K)$ of $K$ is contained in the union of the conormal bundles of the strata of the stratification, the restriction of $H^jK$ to each stratum is not locally constant in general. For example, let $D=\{f:=x^3+x^2z^2-y^2=0\}\subset\C^3$, and consider the stratification of $D$ defined by $S=\{x=y=0\}$, $S'=D\setminus S$. An easy direct verification shows that this stratification is (a)-regular. The singularities of $D$ are resolved by the blow-up $\rho:D'\to D$ along $S$ and there are coordinates $u,v$ of $D'$ such that $\rho^*x=u^2-v^2,\rho^*y=u^3-uv^2,\rho^*z=v$ with $u=\rho^*(y/x)$. Note that $\rho$ is the normalization and $\rho^{-1}(0)$ consists of a point $0'$. Let $K=\rho_*\C_{D'}$. Then $\rank K_x=2$ if and only if $x\in S\setminus\{0\}$. In particular, $K|_S$ is not a local system. On the other hand, the dimension of $CV(K)\cap T_0^*X$ is at most $2$ by the estimation of the characteristic variety of the direct image (see \cite{KS} or \cite{D2}, 4.3.3) because the rank of $d\rho$ at $0'$ is $1$. For $x\in S\setminus\{0\}$, $CV(K)\cap T_x^*X$ consists of two lines spanned respectively by $dg_1$ and $dg_2$ where $g_1,g_2$ are functions defining the local irreducible components of $D$ at $x$. Therefore $CV(K)$ is the closure of the conormal bundle of $S'$ by the involutivity of characteristic varieties (see e.g. \cite{KS}).
\end{rk}

\begin{rk} \label{rk4} \rm
The above example does not give a $\mu$-constant deformation when $z$ is viewed as a parameter, and the $a_f$-condition is not satisfied for $f$. In general, it is known that the converse of a result of \cite{LS} holds, i.e. for a holomorphic function $f$ defined on a neighborhood $X$ of the origin of $\C^n\times\C^r$, the restriction $f_t$ of $f$ to $X\cap\C^n\times\{t\}$ is $\mu$-constant if the stratification of $D=\{f=0\}\subset X$ defined by $S=X\cap\{0\}\times\C^r$ and $S'=D\setminus S$ satisfies the $a_f$-condition (assuming $S'$ smooth). This easily follows from the fact that $f_t$ has a critical point at $x\ne 0$ if and only if the tangent space of $\{f=a\}$ at $(x,t)$ contains $\C^n\times\{0\}$ where $a=f_t(x)$ and $(x,t)\in X$ is sufficiently near $0$. 
\end{rk}

\section{Generalization of Theorem \ref{thm=2} } \label{s4}

\subsection{Deformation to the normal cone} \label{deform}
Let $X$ be a complex manifold, and $Z$ be a closed submanifold. Let $$\hbox{$\X=\Specan_X(\bigoplus_{i\in\Z}I_Z^{-i}\otimes t^i)$},$$ where $I_Z$ is the ideal sheaf of $Z$ in $X$ and $I_Z^{-i}=\OO_{\X}$ for $i\ge 0$.
Note that $\X$ is naturally identified with an open subset of the blow-up of $X\times\C$ along $Z\times\{0\}$. We have the projection $p:\X\to\C$ defined by $t$, and $p^{-1}(\C^*)= X\times\C^*$. Therefore $\X$ gives a deformation of $X$ to the normal cone $N_{Z/X}$ of $Z$ in $X$ (see \cite {Ver}) because $N_{Z/X}$ is isomorphic to $$\hbox{$p^{-1}(0)= \Specan_Z(\bigoplus_{i\le 0}(I_Z^{-i}/I_Z^{-i+1})\otimes t^i)$}.$$

Let $$\hbox{$\A=\bigoplus_{i\in\Z}I_Z^{-i}|_Z\otimes t^i$.}$$ This is identified with a subsheaf of $\OO_{\X}|_Z$ where $Z$ is identified with the zero section of $N_{Z/X}=p^{-1}(0)\subset\X$. Note that the stalk $\A_x$ at each $x\in Z$ is noetherian because there is a surjective ring morphism $\OO_{X,x}[t_0,\dots,t_r]\to\A_x$ sending $t_0$ to $1\otimes t$ and $t_i$ to $x_i\otimes t^{-1}$ if we take local coordinates $x_i$ of $X$ such that $x_i$ for $1\le i\le r$ generate the ideal of $Z$. Moreover, $\A_x$ has a regular sequence consisting of $1\otimes t$, $x_i\otimes t^{-1}$ for $1\le i\le r$, and $x_i\otimes 1$ for $i> r$. Therefore its completion by the maximal ideal generated by these elements is isomorphic to the ring of formal power series over $\C$. This implies that $\OO_{\X,x}$ is flat over $\A_x$ because the completion is faithfully flat over $\OO_{\X,x}$, and is flat over $\A_x$.

\subsection{Specialization} \label{special}
Let $j:X\times\C^*\to\X$ denote the natural inclusion, and $pr:X\times\C^*\to X$ the natural projection. For a regular holonomic $\D_X$-module $M$, consider the open direct image $j_+pr^*M$ of the smooth pull-back $pr^*M$ in the category of regular holonomic $\D$-modules. Note that $j_+pr^*M$ is the localization by $t$ of the pull-back $\rho^*M$ of $M$ by the natural morphism $\rho:\X\to X$, and $$\hbox{$\rho^*M|_Z=(\bigoplus_{i\in\Z}(M\otimes_{\OO_X}I_Z^{-i})|_Z\otimes t^i)\otimes_{\A}\OO_{\X}|_Z$.}$$ Therefore $(j_+pr^*M)|_Z$ is naturally isomorphic to $$\hbox{$(\bigoplus_{i\in\Z}M|_Z\otimes t^i)^{\sim}:=(\bigoplus_{i\in\Z}M|_Z\otimes t^i)\otimes_{\A}\OO_{\X}|_Z$}.$$
Note that $\partial_t$ acts on $(j_+pr^*M)|_Z$, and the kernel of $t\partial_t-i$ is identified with $M|_Z\otimes t^i$, where $\partial_t$ is the vector field on $X\times\C^*$ associated with the coordinate $t$ of $\C^*$ together with the product structure of $X\times\C^*$.

Let $V$ be the filtration of Kashiwara and Malgrange on $M$ along $Z$, and similarly for $j_+pr^*M$ along $p^{-1}(0)$. The specialization of $M$ along $Z$ is defined by $$\psi_tj_+pr^*M.$$
Let $r=\codim_XZ$. It is known (see \cite{BMS} for the algebraic case) that $$\hbox{$V^{\alpha-r+1}(j_+pr^*M)|_Z= (\bigoplus_{i\in\Z}V^{\alpha-i}M|_Z\otimes t^i)^{\sim}$},$$ where $\M^{\sim}$ for an $\A$-module $\M$ in general is defined by $\M\otimes_{\A}\OO_{\X}|_Z$. This isomorphism follows by showing that the filtration defined by the right-hand side satisfies the conditions of the $V$-filtration of Kashiwara and Malgrange. For the proof of this, we can omit the tensor with $\OO_{\X}|_Z$ over $\A$ because this is an exact functor.

Assume that $M$ underlies a mixed Hodge module on $X$. In particular, $M$ has the Hodge filtration $F$. The Hodge filtration on the smooth pull-back $pr^*M$ is given by $pr^*F$, see \cite{S2}. Let $p_0=\min\{p\,|\,F_pM\ne 0\}$. Then $$\hbox{$F_{p_0}(j_+pr^*M)|_Z= (\bigoplus_{i\in\Z}F_{p_0}V^{r-1-i}M|_Z \otimes t^i)^{\sim}$},$$ because $F_{p_0}(j_+pr^*M)=V^0(j_+pr^*M)\cap j_*F_{p_0}(pr^*M)$, see \cite{S1}. Note that $F_{p_0}V^{-i}M=F_{p_0}M$ for $i\gg 0$. Since the tensor with $\OO_{\X}|_Z$ over $\A$ is an exact functor and commutes with intersections of submodules, we get for $\alpha\ge r-1$ $$\hbox{$F_{p_0}V^{\alpha-r+1}(j_+pr^*M)|_Z=(\bigoplus_{i\in\Z}F_{p_0} V^{\alpha-i}M|_Z\otimes t^i)^{\sim}$}.$$ Therefore $F_{p_0}V^{\alpha-i}M|_Z$ is obtained from $F_{p_0}V^{\alpha-r+1}(j_+pr^*M)|_Z$ by restricting to the kernel of $t\partial_t-i$.

\begin{thm} \label{genth2}
Let $D$ be a subvariety of $X$ which is not necessarily reduced nor irreducible. Let $f=(f_1,\dots,f_r)$ be a system of generators of the ideal of $D$, and $i_f:X\to X\times\C^r$ be the graph embedding by $f$.
Assume $T\times\C^r$ is non-characteristic for the specialization of the $\D$-module $i_{f+}\OO_X$ along $X\times\{0\}$, where the normal bundle of $X\times\{0\}$ in $X\times\C^r$ is identified with $X\times\C^r$. Then the  assertion of Theorem \ref{thm=2} holds for $D$.
\end{thm}

\proof

We reduce the assertion to the divisor case applying the above argument to $M=i_{f+}\OO_X$ and $Z=X\times\{0\}\subset X\times\C^r$. Note that the arguments in Section \ref{s3} apply to the case of any mixed Hodge modules including the situation in this section. Note that the vector field $\partial_t$ is defined by using the product structure of $X\times\C$, and induces a vector field on $T\times\C$. Here we may assume that $T$ is defined by a coordinate $z$ of $X$. Then, in order to reduce the assertion to the divisor case as above, it is sufficient to show that the kernel of $t\partial_t-i_0$ and cokernel of the action of $z$ commute on $$\M^{\sim}:=\hbox{$(\bigoplus_{i\in\Z}F_{p_0}V^{\alpha-i}M|_Z\otimes t^i)^{\sim}$.}$$ Indeed, if we take first the cokernel of the action of $z$ we get the intersection of $F$ and $V$ for modules over $T$ by the proof of Theorem \ref{thm=2}.

Let $x\in T\times\{0\}\subset\X$ where $\X$ is an open subvariety of the blow-up of $X\times\C^r\times\C$ along $X\times\{0\}$ and is identified with $X\times\C^r\times\C$.
Using the snake lemma applied to $$\CD0 @>>> \M^{\sim}_x @>{z}>> \M^{\sim}_x @>>> \M^{\sim}_x/z\M^{\sim}_x @>>> 0\\@. @VV{t\partial_t-i_0}V @VV{t\partial_t-i_0}V @VV{t\partial_t-i_0}V\\0 @>>> \M^{\sim}_x @>{z}>> \M^{\sim}_x @>>> \M^{\sim}_x/z\M^{\sim}_x @>>> 0\\
\endCD$$ the above commutativity is reduced to the injectivity of the action of $z$ on the cokernel of $t\partial_t-i_0$ on $\M^{\sim}_x$, and to the canonical isomorphism $$\hbox{Coker}(t\partial_t-i_0:\M^{\sim}_x\to \M^{\sim}_x)=\M^{i_0}_x,$$ where $\M^i=F_{p_0}V^{\alpha-i}M|_Z$. This is further reduced to the surjectivity of the action of $t\partial_t-i_0$ on $$\hbox{$\prod_{i\ne i_0}(\M^i_x \otimes t^i)\bigcap\M^{\sim}_x$}.$$ Here $\M^{\sim}_x$ is viewed as a completion of $\M_x$ in some topology, and is identified with a vector subspace of $\prod_{i}(\M^i_x\otimes t^i)$. This surjectivity is reduced to the case $\M=\A$ and $\M^{\sim}=\OO_{\X}|_Z$ taking homogeneous generators $v_j$ of $\M_x$ (i.e. $v_j\in\M_x^{i_j}\otimes t^{i_j}$).
For $\M=\A$ we take local coordinate system $(x'_1,\dots,x'_n;y'_1,\dots,y'_r;t')$ of $\X$ which is related to a local coordinate system $(x_1,\dots,x_n;y_1,\dots,y_r;t)$ of $X\times\C^r\times\C$ by $x_i=x'_i, y_j= y'_jt', t=t'$. Then the vector field $t\partial_t$ on $\X$ is expressed by $t'\partial_{t'}-\sum_jy'_j\partial_{y'_j}$, and we can prove the surjectivity using the grading of $\C\{x',y',t'\}$ such that $\deg x'_i=0, \deg y'_j=-1, \deg t'=1$.
Therefore the assertion is reduced to the divisor case for arbitrary mixed Hodge modules as above, and follows from the proof of Theorem \ref{thm=2}. This completes the proof of Theorem \ref{genth2}.

\medskip
For the moment, it is not clear whether the assertion of Theorem \ref{genth2} holds under the assumption that $T$ is transversal to any stratum of a Whitney regular stratification of $D$ in the case $D$ is reduced.

\end{document}